\input amstex\documentstyle{amsppt}  
\pagewidth{12.5cm}\pageheight{19cm}\magnification\magstep1
\topmatter
\title Some power series involving involutions in Coxeter groups\endtitle
\author G. Lusztig\endauthor
\address{Department of Mathematics, M.I.T., Cambridge, MA 02139}\endaddress
\thanks{Supported in part by National Science Foundation grant DMS-1303060 and by a Simons Fellowship.}
\endthanks
\endtopmatter   
\document

\define\pos{\text{\rm pos}}

\define\si{\sim}

\define\hM{\hat M}

\define\part{\partial}
\define\emp{\emptyset}

\define\n{\notin}

\define\m{\mapsto}

\define\sub{\subset}    

\define\T{\times}

\define\nl{\newline}
\redefine\i{^{-1}}
\define\fra{\frac}

\define\bst{\bigstar}

\redefine\b{\beta}

\define\g{\gamma}
\redefine\d{\delta}

\define\io{\iota}

\define\ph{\phi}

\define\s{\sigma}

\define\x{\xi}

\define\kk{\bold k}

\define\FF{\bold F}
\define\GG{\bold G}

\define\II{\bold I}

\define\NN{\bold N}

\define\PP{\bold P}
\define\QQ{\bold Q}

\define\ZZ{\bold Z}
\define\XX{\bold X}

\define\ca{\Cal A}
\define\cb{\Cal B}

\define\co{\Cal O}
\define\cp{\Cal P}

\define\car{\Cal R}

\define\cx{\Cal X}

\define\fH{\frak H}

\define\bul{\bullet}

\define\BAR{L1}
\define\HEC{L2}
\define\LV{LV}
\define\MW{MW}

\head Introduction\endhead
\subhead 0.1\endsubhead
Let $(W,S)$ be a Coxeter group such that the canonical set of generators $S$ of $W$ is finite and let 
$w\m w^*$ be an involutive automorphism of $W$ preserving $S$. Let $l:W@>>>\NN$ be the usual length function
and let $\le$ be the standard partial order on $W$. The Poincar\'e series of $W$ is the formal power series
$\PP(u)=\sum_{w\in W}u^{l(w)}\in\ZZ[[u]]$ where $u$ is an indeterminate. We set 
$\PP_*(u)=\sum_{w\in W;w=w^*}u^{l(w)}\in\ZZ[[u]]$. Let $\II_*=\{w\in W;w^*=w\i\}$. It is easy to show (see 
for example \cite{\BAR}) that there is a unique function $\ph:\II_*@>>>\NN$ such that $\ph(1)=0$ and such 
that for any $w\in\II_*$ and any $s\in S$ with $sw<w$ we have $\ph(w)=\ph(sw)+1$ (if $sw=ws^*$) and 
$\ph(w)=\ph(sws^*)$ (if $sw\ne ws^*$). We set 
$$\car(u)=\sum_{z\in\II_*}u^{l(z)}(\fra{u-1}{u+1})^{\ph(z)}\in\ZZ[[u]].$$
(The definition of $\car(u)$ appeared in \cite{\BAR} in the case where $W$ is finite, but the definition
clearly makes sense in the general case.) We say that $W,*$ has property $\XX$ if
$$\car(u)=\PP(u^2)/\PP_*(u).$$
In \cite{\BAR} it was proved (case by case) that $W,*$ has property $\XX$ if $W$ is finite. In \cite{\MW} it
is proved that $W,1$ has property $\XX$ if $W$ is an affine Weyl group of type $A$ and the question of the 
validity of $\XX$ for any $W,*$ is raised. We state the following result.

\proclaim{Theorem 0.2} $W,*$ has property $\XX$.
\endproclaim
The proof is given in 1.11. The most complicated part of the proof is a property stated in Proposition 1.8.
In Section 2 we give an alternative proof of this property assuming that $W$ is a Weyl group. (The same 
proof applies in the case where $W$ is the ``Weyl group'' associated to a Kac-Moody Lie algebra and $*$ is 
induced by an involutive automorphism of that Lie algebra.)

\head 1. The polynomials $X^z_y$\endhead
\subhead 1.1\endsubhead
Let $W,S,*,l,\II_*$ be as in 0.1. Let $u$ be an indeterminate.

Let $\ca=\ZZ[u,u\i]$. Let $\fH$ be the free $\ca$-module with basis 
$(T_w)_{w\in W}$ with the unique $\ca$-algebra structure with unit $T_1$ such that 

$T_wT_{w'}=T_{ww'}$ if $l(ww')=l(w)+l(w')$ and

$(T_s+1)(T_s-u^2)=0$ for all $s\in S$.
\nl
(This is an Iwahori-Hecke algebra.) Let $M$ be the free $\ca$-module with basis $\{a_w;w\in\II_*\}$. 
According to \cite{\LV} (when $W$ is a Weyl group) and \cite{\BAR} (in the general case) there is a unique 
$\fH$-module structure on $M$ (extending the $\ca$-module structure) such that for any $s\in S$ and any 
$z\in\II_*$ we have

(i) $T_sa_z=ua_z+(u+1)a_{sz}$ if $sz=zs^*>z$;

(ii) $T_sa_z=(u^2-u-1)a_z+(u^2-u)a_{sz}$ if $sz=zs^*<z$;

(iii) $T_sa_z=a_{szs^*}$ if $sz\ne zs^*>z$;

(iv) $T_sa_z=(u^2-1)a_z+u^2a_{szs^*}$ if $sz\ne zs^*<z$.

\proclaim{Proposition 1.2} There is a unique collection of polynomials $X^z_y\in\ZZ[u]$ indexed by
$(z,y)\in\II_*\T W$ such that 

($\bst$) $X^z_1=\d_{z,1}$ for all $z\in\II_*$ 
\nl
and such that for any $s\in S,y\in W$ with $y>ys$ and any $z\in\II_*$ we have

(i) $X^z_y=uX^z_{ys}+(u+1)X^{sz}_{ys}$ if $sz=zs^*>z$;

(ii) $X^z_y=(u^2-u-1)X^z_{ys}+(u^2-u)X^{sz}_{ys}$ if $sz=zs^*<z$;

(iii) $X^z_y=X^{szs^*}_{ys}$ if $sz\ne zs^*>z$;

(iv) $X^z_y=(u^2-1)X^z_{ys}+u^2X^{szs^*}_{ys}$ if $sz\ne zs^*<z$.
\nl
Moreover we have $X^z_y\in u^{l(z)}\ZZ[u]$ for any $(z,y)\in\II_*\T W$.
\endproclaim
Let $\x:M@>>>\ca$ be the $\ca$-linear function such that $\x(a_z)=\d_{z,1}$. For $(z,y)\in\II_*\T W$ we set
$X^z_y=\x(T_ya_z)\in\ca$. Then clearly $\bst$ holds.
If $s\in S,y\in W$ satisfy $y>ys$ then $T_y=T_{ys}T_s$; hence applying $T_{ys}$ to
the formulas 1.1(i)-(iv) we see that $X^z_y$ satisfy (i)-(iv) (but they are in $\ca$ 
instead of $\ZZ[u]$). Now from $\bst$ and (i)-(iv) we see by induction on $l(y)$ that $X^z_y$ are uniquely 
determined and that $X^z_y\in u^{l(z)}\ZZ[u]$ for any $(z,y)\in\II_*\T W$. The proposition is proved.

\subhead 1.3\endsubhead
The $\ca$-module $M$ in 1.1 can be viewed as an $\ca$-submodule of $\hM$ which consists of all
formal $\ca$-linear combinations of the elements $\{a_z;z\in\II_*\}$. The $\fH$-module structure on $M$
extends in an obvious way to an $\fH$-module structure on $\hM$. We extend $\x:M@>>>\ca$ to a $\ca$-linear 
map $\hM@>>>\ca$ (denoted again by $\x$) by $\sum_zc_za_z\m c_1$ (here $c_z\in\ca$). We show:
$$T_y(\sum_za_z)=v^{2l(y)}\sum_za_z\text{ for any }y\in W.\tag a$$
It is enough to prove this in the case where $y=s\in S$. For any $z\in\II_*$ we set
$s\bul z=sz$ if $sz=zs^*$ and $s\bul z=szs^*$ if $sz\ne zs^*$. 
From the definition for any $z\in\II_*$ such that $sz<z$ we have
$$T_s(a_z+a_{s\bul z})=u^2(a_z+a_{s\bul z})$$
hence 
$$T_s(\sum_{z\in\II_*}a_z)=\sum_{z\in\II_*;sz<z}T_s(a_z+a_{s\bul z})=
\sum_{z\in\II_*;sz<z}u^2(a_z+a_{s\bul z})=u^2\sum_{z\in\II_*}a_z.$$
This completes the proof of (a).

Applying $\x$ to (a) we obtain for any $y\in W$:
$$\sum_{z\in\II_*}X^z_y=v^{2l(y)}.\tag b$$
Note that in the left hand side we have $X^z_y=0$ for all but finitely many $z$. Indeed 
if $(z,y)\in\II_*\T W$ and $l(z)>2l(y)$, then $X^z_y=0$. (This can be seen from 1.7 or directly from
definitions by induction on $l(y)$.) 

\proclaim{Proposition 1.4} Let $s\in S$, $z\in\II_*$ be such that $sz<z$ and let $y\in W$. We have:

(i) $(u+1)X^z_y=-uX^{sz}_y+X^{sz}_{ys}$ if $sz=zs^*,ys>y$;

(ii) $(u+1)X^z_y=u^2X^{sz}_{ys}+(u^2-u-1)X^{sz}_y$ if $sz=zs^*,ys<y$;

(iii) $X^z_y=X^{szs^*}_{ys}$ if $sz\ne zs^*,ys>y$;

(iv) $X^z_y=u^2X^{szs^*}_{ys}+(u^2-1)X^{szs^*}_y$ if $sz\ne zs^*,ys<y$.
\endproclaim
From 1.2(i) we obtain $X^{sz}_{ys}=uX^{sz}_y+(u+1)X^z_y$ if $sz=zs^*,ys>y$. Clearly this yields (i). In 
1.2(ii) we substitute $X^z_{sy}=-\fra{u}{u+1}X^{sz}_{sy}+\fra{1}{u+1}X^{sz}_y\in\QQ(u)$ (with $ys<y$) which
follows from (i); we obtain
$$\align&X^z_y=(u^2-u-1)(-\fra{u}{u+1}X^{sz}_{sy}+\fra{1}{u+1}X^{sz}_y)+(u^2-u)X^{sz}_{ys}\\&=
\fra{u^2}{u+1}X^{sz}_{ys}+\fra{u^2-u-1}{u+1}X^{sz}_y\endalign$$
so that (ii) holds. Now (iii) clearly follows from 1.2(iii). In 1.2(iv) we substitute $X^z_{ys}=X^{szs^*}_y$
which follows from (iii); we obtain (iv). The proposition is proved.

From the proposition we deduce the following result.
\proclaim{Corollary 1.5} Let $z\in\II_*$, $y\in W$, $s\in S$. We have:

(i) $X^z_y+X^z_{ys}=u^2(X^{szs^*}_y+X^{szs^*}_{ys})$ if $sz\ne zs^*<z$;

(ii) $(u+1)(X^z_y+X^z_{ys})=(u^2-u)(X^{sz}_y+X^{sz}_{ys})$ if $sz=zs^*<z$.
\endproclaim

\subhead 1.6\endsubhead
For any $m\in M$ we write $m=\sum_{z\in\II_*}(a_z:m)a_z$ where $(a_z:m)\in\ca$ are zero for all but finitely
many $z$. For any $h\in\fH$ we write $h=\sum_{w\in W}[T_w:h]T_w$ where $[T_w:h]\in\ca$ are zero for all but 
finitely many $w$. 

For $z'\in\II_*,s\in S$ we show:

(i) $T_sT_{z'}T_{s^*}=cT_{z'}+c'T_{sz'}$ where $c,c'\in\ZZ[u]-\{0\}$, if $sz'=z's^*$;

(ii) $T_sT_{z'}T_{s^*}=T_{sz's^*}$ if $sz'\ne z's^*>z$;

(iii) $T_sT_{z'}T_{s^*}=u^4T_{sz's^*}+(u^2-1)^2T_{z'}+u^2(u^2-1)T_{sz'}+u^2(u^2-1)T_{z's^*}$ if 
$sz'\ne z's^*<z'$.
\nl
In (i) assume first that $sz'>z'$. Then $T_sT_{z'}T_{s^*}=T_sT_{z's^*}=u^2T_{z'}+(u^2-1)T_{sz'}$. 
Now assume that in (i) we have $sz<z$. Then
$$\align&T_sT_{z'}T_{s^*}=T_sT_sT_{sz'}T_{s^*}=T_sT_sT_{z'}=u^2T_{z'}+(u^2-1)T_sT_{z'}\\&
=u^2T_{z'}+u^2(u^2-1)T_{sz'}+(u^2-1)^2T_{z'}.\endalign$$
Now (ii) is obvious. In (iii) we have
$$\align&T_sT_{z'}T_{s^*}=T_sT_sT_{sz's^*}T_{s^*}T_{s^*}\\&
=u^4T_{sz's^*}+(u^2-1)^2T_{z'}+u^2(u^2-1)T_{sz'}+u^2(u^2-1)T_{z's^*},\endalign$$
as desired.

For $s\in S$ and $z,z'\in\II_*$ we show:

(a) If $(a_{z}:T_sa_{z'})\ne0$ then $[T_z:T_sT_{z'}T_{s^*}]\ne0$.
\nl
If $s'z=z's^*$ then by 1.2 we have $z=z'$ or $z=sz'$ and (a) follows from (i). 

If $sz'\ne z's^*>z'$ then by 1.2 we have $z=sz's^*$ and (a) follows from (ii). 

If $sz'\ne z's^*<z'$ then by 1.2 we have $z=z'$ or $z=sz's^*$ and (a) follows from (iii). 

This proves (a). 

We shall also need a variant of (a). Assume that $s\ne s'$ in $S$ are such that $s'=s^*$, $s,s'$ generate a 
finite (dihedral) subgroup of $W$ and $\s$ is the longest element in that subgroup. Then

(b) If $z,z'\in\II_*$ and  $(a_{z}:T_\s a_{z'})\ne0$ then $[T_z:T_\s T_{z'}T_{\s^*}]\ne0$;

(c) $(a_1:T_\s a_1)=u^{l(\s)}$.
\nl
This can be proved using the results in \cite{\BAR} which describe the action of $T_\s$ on $M$.

For $(z,y)\in\II_*\T W$ we set $c^z_y=(a_z:T_ya_1)$. We have the following result.

\proclaim{Proposition 1.7}If $z\in\II_*,y\in W$ and $c^z_y\ne0$ then $[T_z:T_yT_{y^{*-1}}]\ne0$.
\endproclaim
We argue by induction on $l(y)$. If $y=1$ the result is obvious. Assume now that $y\ne1$. We can find
$s\in S$ such that $sy<y$. We have 
$$c^z_y=\sum_{z'\in\II^*}(a_z:T_sa_{z'})c^{z'}_{sy}.$$
By assumption there exists $z'\in\II^*$ such that $(a_z:T_sa_{z'})\ne0$ and $c^{z'}_{sy}\ne0$. By 1.6(a) we 
then have $[T_z:T_sT_{z'}T_{s^*}]\ne0$ and by the induction hypothesis we have 
$[T_{z'}:T_{sy}T_{y^{*-1}s^*}]\ne0$. We have
$$[T_z:T_yT_{y^{*-1}}]=\sum_{w\in W}[T_z:T_sT_wT_{s^*}][T_w:T_{sy}T_{y^{*-1}s^*}].$$
This is a sum of terms in $\NN[u-1]$ (the set of polynomials in $u-1$ with coefficients in $\NN$) 
and the term corresponding to $w=z'$ is nonzero; hence the sum is nonzero. The proposition is proved.

\proclaim{Proposition 1.8}For any $y\in W$ we have $X^1_y=\d_{y,y^*}u^{l(y)}$.
\endproclaim
From \cite{\HEC, 10.4(a)} we have
$$[T_1:T_{y_1}T_{y_1^{*-1}}]=\d_{y_1,y_1^*}u^{2l(y_1)}\text{ for any }y_1\in W.\tag a$$
We prove the proposition by induction on $l(y)$. If $y=1$ the result is obvious. Assume now that $y\ne1$. 
We can write $y=\s y'$ where $\s,y'\in W$, $l(y)=l(\s)+l(y')$ and the following holds: if $y\ne y^*$ then 
$\s\in S$; if $y=y^*$ then there exists $s,s'$ in $S$ such that $s'=s^*$, $s,s'$ generate a finite subgroup 
of $W$ and $\s$ is the longest element in that subgroup, so that $\s=\s^*$. (We use \cite{\HEC, A1(a)}.) By 
the induction hypothesis we have
$$T_{y'}a_1=\d_{y',y'{}^*}u^{l(y')}a_1+\sum_{z\in\II_*,z\ne1}c^z_{y'}a_z.$$
Since $T_y=T_\s T_{y'}$, it follows (using also (a)) that
$$T_ya_1=\d_{y',y'{}^*}u^{l(y')}T_\s a_1+\sum_{z\in\II_*,z\ne1}c^z_{y'}T_\s a_z,\tag b$$
$$T_yT_{y^{*-1}}=\d_{y',y'{}^*}u^{2l(y')}T_\s T_{\s^*}+\sum_{w\in W;w\ne1}[T_w:T_{y'}T_{y'{}^{*-1}}]
T_\s T_wT_{\s^*}.$$
Note that by (a) we have
$$[T_1:T_\s T_{\s^*}]=\d_{\s,\s^*}u^{2l(\s)}.$$
Hence
$\g:=[T_1:\d_{y',y'{}^*}u^{2l(y')}T_\s T_{\s^*}]$ is $\d_{y',y'{}^*}u^{2l(y')}\d_{\s,\s^*}u^{2l(\s)}$, that 
is $\d_{y,y^*}u^{2l(y)}$. (If $y\ne y^*$ we have $\g=0$ since we have either $\d_{y',y'{}^*}=0$ or 
$\d_{\s,\s^*}=0$. If $y=y^*$ then $\s=\s^*$ hence $y'=y'{}^*$ and $\g=u^{2l(y)}$.)

Assume that for some $z\in\II_*-\{1\}$, we have $c^z_{y'}\ne0$ and $(a_1:T_\s a_z)\ne0$. Then, by 1.7
we have $[T_z:T_{y'}T_{y'{}^*}]\ne0$. Moreover we have $[T_1:T_\s T_zT_{\s^*}]\ne0$. (When $\s\in S$ this
follows from 1.6(a). When $s\n S$ this follows from 1.6(b).)
Thus, $\pi:=[T_1:\sum_{w\in W;w\ne1}d^w_{y',y'{}^{*-1}}T_\s T_wT_{\s^*}]$ is a sum of terms in $\NN[u-1]$, 
at least one of which is $\ne0$, so that $\pi\in\NN[u-1]-\{0\}$. Thus 
$[T_1:T_yT_{y^{*-1}}]=\d_{y,y^*}u^{2l(y)}+\pi$. By (a), $[T_1:T_yT_{y^{*-1}}]=\d_{y,y^*}u^{2l(y)}$. Thus,
$\pi=0$, a contradiction.

This contradiction shows that for any $z\in\II_*-\{1\}$ we have either $c^z_{y'}=0$ or $(a_1:T_\s a_z)=0$.
Hence from (b) we deduce that $(a_1:T_ya_1)=(a_1:\d_{y',y'{}^*}u^{l(y')}T_\s a_1)$. We have 
$$(a_1:T_\s a_1)=\d_{\s,\s^*}u^{l(\s)}.$$
(When $\s\in S$ this follows from 1.1(i)-(iv). When $\s\n S$ this follows from 1.6(c).) We deduce
$$(a_1:T_ya_1)=\d_{y',y'{}^*}u^{l(y')}\d_{\s,\s^*}u^{l(\s)}=\d_{y',y'{}^*}u^{l(y')}\d_{\s,\s^*}u^{l(y)}.$$
It remains to show that 
$$\d_{y',y'{}^*}\d_{\s,\s^*}=\d_{y,y^*}.\tag c$$
If $y\ne y^*$ then we have either $y'\ne y'{}^*$ or $\s\ne\s^*$ hence both sides of (c) are zero. 
If $y=y^*$ then $\s=\s^*$ hence $y'=y'{}^*$ hence both sides of (c) are $1$. Thus (c) holds.
The proposition is proved.

\proclaim{Proposition 1.9}
We have $X^z_y\in u^{l(y)}\ZZ[u]$ for any $(z,y)\in\II_*\T W$.
\endproclaim
This follows by induction on $l(z)$ using 1.4(i)-(iv); to start the induction, we assume that $z=1$ in which
case the result follows from 1.8.

\subhead 1.10\endsubhead
For any $z\in\II_*$ we set
$$X^z=\sum_{y\in W}X^z_y\in\ZZ[[u]].\tag b$$
Note that, by 1.9, the sum in the right hand side converges in $\ZZ[[u]]$. For $z\in\II_*$ we show:
$$X^z=\PP_*(u)u^{l(z)}(\fra{u-1}{u+1})^{\ph(z)}.\tag c$$
We argue by induction on $l(z)$. If $z=1$ we have, using 1.8:
$$X^1=\sum_{y\in W;y=y^*}u^{l(y)}=\PP_*(u).$$
Assume now that $z\ne1$. We can find $s\in S$ such that $sz<z$. From 1.5 we deduce that $X^z=u^2X^{szs^*}$ 
if $sz\ne zs^*$ and $(u+1)X^z=(u^2-u)X^{sz}$ if $sz=zs^*$. Using the induction hypothesis we see that, if 
$sz\ne zs^*$ we have 
$$X^z=\PP_*(u)u^{l(szs^*)+2}(\fra{u-1}{u+1})^{\ph(szs^*)},$$
while if $sz=zs^*$ we have 
$$X^z=\PP_*(u)u^{l(sz)+1}(\fra{u-1}{u+1})^{\ph(sz)+1};$$
the desired result follows.

\subhead 1.11\endsubhead
We prove Theorem 0.2. The sum $\sum_{(z,y)\in\II_*\T W}X^z_y$ is convergent in $\ZZ[[u]]$ since 
$X^z_y\in u^{\max(l(z),l(y))}\ZZ[u]$ (see Propositions 1.2 and 1.9). We can compute this sum in two 
different ways and we get the same result. Thus we have
$$\sum_{z\in\II_*}(\sum_{y\in W}X^z_y)=\sum_{y\in W}(\sum_{z\in\II_*}X^z_y).$$
By 1.3(b), the right hand side is equal to $\sum_{y\in W}u^{2l(y)}=\PP(u^2)$. By 1.4(c), the left hand side 
is equal to
$$\sum_{z\in\II_*}\PP_*(u)u^{l(z)}(\fra{u-1}{u+1})^{\ph(z)}=\PP_*(u)\car(u).$$
Thus we have 
$$\PP(u^2)=\PP_*(u)\car(u).$$
Theorem 0.2 is proved.

\head 2. The case of Weyl groups\endhead
\subhead 2.1\endsubhead
In this section we assume that $W$ is the Weyl group of a connected adjoint simple algebraic group $\GG$
defined and split over the finite field $\FF_p$ with $p$ elements ($p$ is a prime number); we identify $\GG$
with $\GG(\kk)$ where $\kk$ is an algebraic closure of $\FF_p$.  
Let $F:\GG@>>>\GG$ be the ``Frobenius map''; it is an abstract group 
isomorphism whose fixed point set is the group of $\FF_p$-rational points of $\GG$. Let $\cb$ be the 
variety of Borel subgroups of $\GG$. The set of $\GG$-orbits on $\cb\T\cb$ (for the 
simultaneous conjugation action) are naturally indexed by $W$; we denote by $\co_w$ the $\GG$-orbit indexed 
by $w\in W$. The length function $l:W@>>>\NN$ has the 
property that for any $w\in W$ and any $C\in\cb$, the algebraic variety
$\{B\in\cb;(C,B)\in\co_w\}$ is an affine space over $\kk$ of dimension $l(w)$. The subset $S$ of $W$ is then
given by $\{s\in W;l(w)=1\}$. If $B\in\cb$ then $F(B)\in\cb$; moreover, $B\m F(B)$ is a bijection 
$F:\cb@>>>\cb$.

For a finite set $X$ we denote by $|X|$ the cardinal of $X$. If $X\sub X'$ are sets and 
$f:X'@>>>X'$ is a map such that $f(X)\sub X$ we set $X^f=\{x\in X;f(x)=x\}$.

\subhead 2.2\endsubhead
Let $s\in S$. If $B,B'\in\cb$ we write $B\si_sB'$ if $(B,B')\in\co_1\cup\co_s$. This is an equivalence 
relation on $\cb$. A subgroup $P$ of $\GG$ is said to be parabolic of type $s$ if it is the union of all $B$
in a fixed equivalence class for $\si_s$. Let $\cp_s$ be the set of parabolic subgroups of type $s$ of $\GG$.
For $P\in\cp_s$ we set $\cb_P=\{B\in\cb;B\sub P\}$; this is a projective line over $\kk$.

If $P\in\cp_s$ and $B\in\cb$, then there is a unique element $y=\pos(B,P)\in W$ such that $y<ys$ and 
$(B,B')\in\co_y\cup\co_{ys}$ for any $B'\in\cb_P$; moreover we have $(B,B')\in\co_y$ for a unique 
$B'\in\cb_P$, denoted by $P^B$. 

\subhead 2.3\endsubhead
Let $s,s'\in S$ and let $P\in\cp_s,P'\in\cp_{s'}$. There is a unique element $y=\pos(P,P')\in W$ such that 
$y<sy,y<ys'$ and $(B,B')\in\co_y\cup\co_{sy}\cup\co_{ys'}\cup\co_{sys'}$ for any $B\in\cb_P,B'\in\cb_{P'}$.
Assuming in addition that $sy=ys'$, the map $\cb_P@>>>\cb_{P'}$ given by $B\m P'{}^B$ is an isomorphism of 
projective lines with inverse $B'\m P^{B'}$.

\subhead 2.4\endsubhead
We shall fix $C\in\cb^F$. We assume that $*$ has the following property: we can find an involutive 
automorphism of algebraic groups $\io:\GG@>>>\GG$ such that:

$\io$ commutes with $F:\GG@>>>\GG$;

$\io(C)=C$; 

if $(B,B')\in\co_w$ then $(\io(B),\io(B'))\in\co_{w^*}$.
\nl
Note that $\io$ defines an involution $\cb@>>>\cb$ of algebraic varieties (denoted again by $\io$). It
commutes with $F:\cb@>>>\cb$.

\subhead 2.5\endsubhead
For $(z,y)\in\II_*\T W$ we set 
$$\cx^z_y=\{B\in\cb^{F^2};(B,F(\io(B)))\in\co_z,(C,B)\in\co_y\}.$$
This is a finite set.

\proclaim{Proposition 2.6} Let $(z,y)\in\II_*\T W$, $s\in S$. Assume that $sz<z$. We have

(i) $(p+1)|\cx^z_y|=-p|\cx^{sz}_y|+|\cx^{sz}_{ys}|$ if $sz=zs^*$, $y<ys$;

(ii) $(p+1)|\cx^z_y|=p^2|\cx^{sz}_{ys}|+(p^2-p-1)|\cx^{sz}_y|$ if $sz=zs^*$, $y>ys$;

(iii) $|\cx^z_y|=|\cx^{szs^*}_{ys}|$ if $sz\ne zs^*$, $y<ys$;

(iv) $|\cx^z_y|=p^2|\cx^{szs^*}_{ys}|+(p^2-1)|\cx^{szs^*}_y|$ if $sz\ne zs^*$, $y>ys$.
\endproclaim
In the setup of (iii), (iv) we have $l(z)=l(s)+l(szs^*)+l(s^*)$ hence for any $B\in\cx^z_y$ we have
$(B,\b_B)\in\co_s$, $(\b_B,\b'_B)\in\co_{szs^*}$, $(\b'_B,F(\io(B)))\in\co_{s^*}$ for uniquely defined 
$\b_B,\b'_B$ in $\cb$. It follows that $(B,F(\io(\b'_B)))\in\co_s$, $(F(\io(\b'_B)),F(\io(\b_B)))\in
\co_{szs^*}$, $(F(\io(\b_B)),B)\in\co_{s^*}$ and from the uniqueness we see that $\b'_B=F(\io(\b_B))$, 
$\b_B=F(\io(\b'_B))$; hence $F^2(\b_B)=\b_B$.

Assume that we are in the setup of (iii). We have a bijection $\cx^z_y@>>>\cx^{szs^*}_{ys}$ given by 
$B\m\b_B$. Hence (iii) holds.

Assume that we are in the setup of (iv). We have a map $\cx^z_y@>>>\cx^{szs^*}_y\cup\cx^{szs^*}_{ys}$ given 
by $B\m\b_B$. Its fibre over a point in $\cx^{szs^*}_y$ has cardinal $p^2-1$ while its fibre over a point in
$\cx^{szs^*}_{ys}$ has cardinal $p^2$. We deduce that (iv) holds.

In the rest of the proof we assume that $sz=zs^*$. We set $y'=y$ if $y<ys$ and $y'=ys$ if $ys<y$. Let 
$$Y=\{P\in\cp_s;F^2(P)=P,\pos(P,F(\io(P)))=sz,\pos(C,P)=y'\}.$$
For $P\in Y$ we have $P^C\in\cb_P$, $(F(P))^C=F(P^C)\in\cb_{F(P)}$. 

We define a map $\cx^z_{y'}\cup\cx^{sz}_{y'}\cup\cx^z_{y's}\cup\cx^{sz}_{y's}@>>>Y$ by $B\m P$ where 
$B\sub P\in\cp_s$. The fibre of this map over $P\in Y$ is $\{\b\in\cb^{F^2};B\sub P\}$ hence it has cardinal
$p^2+1$.

We define a bijection $\cx^z_{y'}\cup\cx^{sz}_{y'}@>>>Y$ by $B\m P$ where $B\sub P\in\cp_s$. 

We define a map $\cx^{sz}_{y'}\cup\cx^{sz}_{y's}@>>>Y$ by $B\m P$ where $B\sub P\in\cp_s$. The fibre of this
map over $P\in Y$ is $\{\b\in\cb^{F^2};\b\sub P,(\b,F(\b))\in\co_{sz}\}$ which can be viewed as a set of 
$\FF_p$-rational points of the projective line $\cb_P$ with Frobenius map $\b\m P^{F(\b)}$; hence it has 
cardinal $p+1$.

We see that

$|\cx^z_{y'}|+|\cx^{sz}_{y'}|+|\cx^z_{y's}|+|\cx^{sz}_{y's}|=(p^2+1)|Y|$,

$|\cx^z_{y'}|+|\cx^{sz}_{y'}|=|Y|$,

$|\cx^{sz}_{y'}|+|\cx^{sz}_{y's}|=(p+1)|Y|$.
\nl
Thus we have 

$|\cx^x_{y'}|+|\cx^{sz}_{y'}|+|\cx^z_{y's}|+|\cx^{sz}_{y's}|=(p^2+1)(|\cx^z_{y'}|+|\cx^{sz}_{y'}|)$,

$|\cx^{sz}_{y'}|+|\cx^{sz}_{y's}|=(p+1)(|\cx^z_{y'}|+|\cx^{sz}_{y'}|)$.
\nl
This implies that (i),(ii) hold.

\proclaim{Proposition 2.7} Let $(z,y)\in\II_*\T W$, $s\in S$. Assume that $ys<y$. We have

(i) $\cx^z_y=p\cx^z_{ys}+(p+1)\cx^{sz}_{ys}$ if $sz=zs^*>z$;

(ii) $\cx^z_y=(p^2-p-1)\cx^z_{ys}+(p^2-p)\cx^{sz}_{ys}$ if $sz=zs^*<z$;

(iii) $\cx^z_y=\cx^{szs^*}_{ys}$ if $sz\ne zs^*>z$;

(iv) $\cx^z_y=(p^2-1)\cx^z_{ys}+p^2\cx^{szs^*}_{ys}$ if $sz\ne zs^*<z$.
\endproclaim
The proposition is deduced from 2.6 in the same way as 1.4 was deduced from 1.2 (or rather in reverse).

\proclaim{Proposition 2.8} For any $z\in\II_*$ we have $|\cx^z_1|=\d_{z,1}$.
\endproclaim
We have
$$\cx^z_1=\{B\in\cb^{F^2};(B,F(\io(B)))\in\co_z,C=B\}.$$
Since $C=F(\io(C))$ we have $(C,F(\io(C)))\in\co_1$ hence $\cx^z_1$ is a single point $C$ if $z=1$ and is 
empty if $z\ne1$. The proposition is proved.

\proclaim{Proposition 2.9} For any $y\in W$ we have $|\cx^1_y|=\d_{y,y^*}p^{l(y)}$.
\endproclaim
We have
$$\align&\cx^1_y=\{B\in\cb^{F^2};B=F(\io(B)),(C,B)\in\co_y\}\\&=\{B\in\cb;B=F(\io(B)),(C,B)\in\co_y\}.
\endalign$$
If $B\in\cx^1_y$ then $(F(\io(C)),F(\io(B))\in\co_{y^*}$
hence $(C,B)\in\co_{y^*}$ so that $y=y^*$. Thus, if $\cx^1_y\ne\emp$ then $y=y^*$.
Assume now that $y=y^*$. Consider the $l(y)$-dimensional affine space $\{B\in\cb;(C,B)\in\co_y\}$ over $\kk$.
Since $y=y^*$, this affine space is stable under the map $B\m F(\io(B))$ which is the Frobenius map for an
$\FF_p$-rational structure whose fixed point set is exactly $\cx^1_y$. It follows that $|\cx^1_y|=p^{l(y)}$.
The proposition is proved.

\proclaim{Proposition 2.10} For any $(z,y)\in\II_*\T W$ we have $|\cx^z_y|=X^z_y|_{u=p}$.
\endproclaim
By 1.2 and 2.8, the two sides of the equality in the proposition satisfy the same inductive formulas; they
also satisfy the same initial conditions $|\cx^z_1|=X^z_1|_{u=p}=\d_{z,1}$ (see 2.8 and 1.2). Hence they are
equal. The proposition is proved.

\proclaim{Proposition 2.11} For any $y\in W$ we have $X^1_y=\d_{y,y^*}u^{l(y)}$.
\endproclaim
Combining Propositions 2.9 and 2.10 we see that $X^1_y|_{u=p}=\d_{y,y^*}u^{l(y)}|_{u=p}$. Since two 
polynomials in $u$ which take equal values at any prime number are equal, we deduce that
$X^1_y=\d_{y,y^*}u^{l(y)}$. The proposition is proved.

Note that this proof gives an alternative approach to Proposition 1.8 in our case.

\enddocument